\documentclass[12pt]{amsart}

\usepackage{amscd}
\usepackage{amssymb}

\allowdisplaybreaks[1]

\newtheorem{thm}{Theorem}[section]
\newtheorem{lem}[thm]{Lemma}
\newtheorem{cor}[thm]{Corollary}

\theoremstyle{definition}

\newtheorem{rem}[thm]{Remark}
\newtheorem{defn}[thm]{Definition}

\theoremstyle{remark}

\numberwithin{equation}{section}

\def\N{{\mathbb N}}
\def\C{{\mathbb C}}
\def\Ma{{\mathbb M}}
\def\R{{\mathbb R}}

\def\M{\text{$(M,\Delta)$}} 
\def\J{\text{$J_\varphi$}}
\def\J^{\text{$\hat{J}$}}

\DeclareMathOperator{\spa}{span}
\DeclareMathOperator{\Tr}{Tr}

\DeclareMathOperator{\id}{id}
\DeclareMathOperator{\re}{Re}

\def\al{\alpha}

\def\la{\lambda}

\def\vph{\varphi}
\def\om{\omega}
\def\si{\sigma}
\def\Si{\Sigma}

\def\De{\Delta}
\def\vDe{\varDelta}

\def\vLa{\varLambda}

\def\oti{\otimes}

\def\lan{\langle}
\def\ran{\rangle}
\def\haA{\hat{A}}
\def\haM{\hat{M}}
\def\haD{\hat{\Delta}}

\def\haJ{\hat{J}}

\def\haW{\widehat{W}}
\def\havLa{\hat{\vLa}}
\def\havph{\hat{\vph}}

\allowdisplaybreaks[4]

\begin{document}

\title[Amenable Discrete Quantum Groups]
{Amenable Discrete Quantum Groups}

\author[Reiji Tomatsu]{Reiji Tomatsu}
\address{Department of Mathematical Sciences,
University of Tokyo, 3-8-1 Komaba, Meguro, Tokyo 153-8914, Japan}
\email{tomatsu\char`\@ms.u-tokyo.ac.jp}


\begin{abstract}
Z.-J. Ruan has shown that several amenability conditions are 
all equivalent in the case of discrete Kac algebras. 
In this paper, we extend this work to the case of 
discrete quantum groups. 
That is, we show that a discrete quantum group, where we do not assume 
its unimodularity, has an invariant mean 
if and only if it is strongly Voiculescu amenable. 
\end{abstract}

\maketitle

\section{Introduction}
In this paper, we study amenability of non-Kac type 
discrete quantum groups. 
We use notions of discrete quantum groups or its dual compact 
quantum groups introduced in \cite{Effros Ruan}, \cite{Woronowicz2}, 
for example. 
Amenability is defined as a generalization of the group case, 
that is, by the existence of an invariant mean. 
In a discrete group case, it is known that amenability 
is characterized by several conditions 
(see \cite{Kye1} for its survey). 
The first step to its generalization for discrete quantum groups has been 
made by Z.-J. Ruan \cite{Ruan1} 
under the tracial condition of the Haar weight, 
that is, the Kac algebra condition. 
In particular, he has shown that amenability is equivalent to 
strong Voiculescu amenability. 
For general quantum groups, a generalization 
of Ruan's theorem has been investigated by 
E. B\'{e}dos, R. Conti, G.J. Murphy and L. Tuset 
in \cite{Bedos Conti Tuset}, \cite{Bedos Murphy Tuset1}, 
\cite{Bedos Murphy Tuset2}, \cite{Bedos Murphy Tuset3}. 
They have shown that strong Voiculescu amenability implies 
amenability. 
Our main theorem (Theorem \ref{main}) 
says that both the notions are equivalent for 
general discrete quantum groups. 
We should mention that all the implications except for the above 
converse one in Theorem \ref{main} 
have been already known in the pioneering works 
\cite{Bedos Conti Tuset}, \cite{Bedos Murphy Tuset1}, 
\cite{Bedos Murphy Tuset2}, \cite{Bedos Murphy Tuset3}, 
however we give a proof in order to give another proof of 
nuclearity of dual compact quantum groups. 
After this work was done, we learned from S. Vaes that 
E. Blanchard and he also have proved equivalence between 
amenability and strong Voiculescu amenability. 

\textbf{Acknowledgements}.
The author is highly grateful to his supervisor Yasuyuki Kawahigashi, 
Yoshiomi Nakagami and Stefaan Vaes for all the discussions and encouragement.

\section{Notations for quantum groups}
On symbols of tensor products (minimal tensor product or tensor product 
von Neumann algebra), the same notation $\otimes$ 
is used throughout this paper. 
A flip unitary on tensor product Hilbert spaces $H\oti H$ 
is denoted by $\Si$. 
For an weight $\theta$ on a von Neumann algebra $N$, 
define the left ideal $n_\theta$ 
and $*$-subalgebra $m_\theta$ by 
$n_\theta=\{x\in N \mid \theta(x^*x)<\infty\}$ and 
$m_\theta=n_\theta^* n_\theta$. 
$m_\theta^+$ means $m_\theta\cap N_+$. 

\subsection{locally compact quantum groups and its duals} 

We adopt the definition of locally compact quantum groups 
advocated in \cite{Kustermans Vaes} as follows. 

\begin{defn}
A pair $(M,\Delta)$ is called a (von Neumann algebraic) 
locally compact quantum group 
when it satisfies the following conditions.
\begin{enumerate}
  
\item 
$M$ is a von Neumann algebra and 
$\Delta:M\longrightarrow M\otimes M$ is a unital normal 
$*$-homomorphism satisfying the coassociativity relation: 
$(\Delta\otimes\iota)\Delta =(\iota\otimes\Delta)\Delta$.

\item 
There exist two faithful normal semifinite weights $\varphi$ 
and $\psi$ 
which satisfy 
$\varphi((\omega\otimes\iota)\Delta(x))
=\omega(1)\varphi(x)$ 
for all $x\in m_{\vph}^{+}$, $\omega\in M_{*}^{+}$ 
and 
$\psi((\iota\otimes\omega)\Delta(x))=\omega(1)\psi(x)$ 
for all $x\in m_{\psi}^{+}$, $\omega\in M_{*}^{+}$.     
\end{enumerate}
\end{defn}

The von Neumann algebra $M$ is realized in $B(H)$ via the GNS 
representation associated to $\varphi$, $\{H,\vLa\}$ 
where $\vLa$ is a map from the left ideal 
$n_\vph=\{x\in M\mid \vph(x^*x)<\infty\}$ to the Hilbert space $H$. 
On the tensor product $H\oti H$, the multiplicative unitary $W$ 
is defined as follows, for $x,y\in n_\vph$, 
\[
W^* (\vLa (x) \otimes \vLa (y)) 
=(\vLa \otimes \vLa )
(\varDelta (y)(x\otimes 1)).
\]
It satisfies the pentagonal equality 
$W_{12}W_{13}W_{23}=W_{23}W_{12}$. 
We often use a $C^*$-subalgebra $A$ of $M$ which is defined as 
the norm closure of the linear space 
$\{(\id\oti\om)(W)\mid \om\in B(H)_*\}$. 
It will be considered a continuous function part of $M$. 

This unitary $W$ also plays a role in defining 
the dual locally compact quantum group $(\haM,\haD)$. 
The von Neumann algebra $\haM$ is 
the $\sigma$-weak closure of 
$\{\la(\om) \mid  \omega \in B(H)_* \}$, 
where $\la(\om)=(\omega\oti\id)(W)$. 
Set $\haW=\Si W^* \Si$ and its coproduct is given 
by $\haD(x)=\haW^*(1\oti x)\haW$. 
As above, a $C^*$-subalgebra $\haA$ of $\haM$ is also defined by 
the norm closure of the linear space 
$\{(\id\oti\om)(\haW)\mid \om\in B(H)_*\}$. 

The left invariant weight $\havph$ on $\haM$ 
is characterized by the following property. 
For $\om\in M_*$, if a vector $\xi$ meets 
$\om(x^*)=\lan \xi \mid \vLa(x) \ran$ for all $x\in n_\vph$, 
then $\havph(\la(\om)^*\la(\om))=||\xi||^2$. 
Note that $H$ also becomes a GNS Hilbert space for $\havph$ via 
$\havLa(\la(\om))=\xi$. 
Denote the set of such $\om$ by $\mathcal{I}$, 
that is a dense subspace of $M_*$. 

Let $J$ and $\haJ$ be the modular conjugation for $\vph$ 
and $\havph$. 
Then the following useful equalities hold, 
\[
(\haJ\oti J)W(\haJ\oti J)=W^*,
\quad
(J\oti \haJ)\haW(J\oti \haJ)=\haW^*.
\]
The modular operator and the modular automorphism of $\vph$ 
is denoted by $\vDe_\vph$ and $\si^\vph$, respectively. 
The autopolar of $\vph$, $\mathcal{P}_\vph^\natural$ is defined by 
the norm closure of $\{xJ\vLa(x)\mid x\in n_\vph \cap n_\vph^*\}$. 
Then any normal state on $M$ is of the form $\om_\xi$ with 
$\xi\in \mathcal{P}_\vph^\natural$ and such a vector 
is unique (\cite{Haagerup}). 

\subsection{Discrete and compact quantum groups}

A locally compact quantum group $(M,\De)$ is called 
\textit{discrete} if $\havph(1)<\infty$. 
Then its dual $(\haM,\haD)$ is called \textit{compact}. 
In this case, the state condition $\havph(1)=1$ is always assumed. 
$\havph$ has the left and right invariance. 
About them, the basic references are \cite{Effros Ruan}, 
\cite{Van Daele} and \cite{Woronowicz2}. 
Then $M$ becomes a direct sum von Neumann algebra 
of matrix algebras, say 
$M=\oplus_{\al\in I} \Ma_{n_\al}(\C)$. 
Then its left invariant weight $\vph$ is decomposed as 
$\oplus_\al \Tr_\al h_\al$, 
where $\Tr_\al$ are usual non-normalized trace. 
For $\al$, fix a matrix unit $\{e(\al)_{i,j}\}_{i,j\in I}$ 
which diagonalize $h_\al$ as 
$h_\al=\sum_{i\in I}\nu(\al)_i\, e(\al)_{i,i}$. 
The positive affiliated operator 
$h'=\sum_{\al\in I}\Tr_\al(h_\al)^{-1}h_\al$ is group-like, i.e. 
$\vDe(h')=h'\oti h'$. 
Hence for the modular automorphism $\si^\vph$ and the coproduct 
$\De$, we have for $t\in \R$, 
\[
\De\circ \si_t^\vph=(\si_t^\vph\oti \si_t^\vph )\circ \De. 
\]
Note that the above equality does not hold for general cases. 

\section{Amenability of quantum groups}
We begin with the following well-known definition. 

\begin{defn}
Let $(M,\Delta)$ be a locally compact quantum  group.
\begin{enumerate}
  
\item 
A state $m$ of $M$ is called a left invariant mean 
if 
$m((\omega\otimes\iota)(\Delta(x)))=\omega(1)m(x)$ 
for all $\omega \in M_*$ and $x\in M$. 

\item 
A state $m$ of $M$ is called a right invariant mean 
if $m((\iota\otimes\omega)(\Delta(x)))=\omega(1)m(x)$ 
for all $\omega \in M_*$ and $x\in M$.
  
\item A state $m$ of $M$ is called an invariant mean 
if $m$ is a left and right invariant mean.

\end{enumerate}
\end{defn}

\begin{rem}
If $(M,\Delta)$ has a left invariant mean, 
it also has an invariant mean 
(see \cite[Proposition 3]{Desmedt Quaegebeur Vaes} for its proof). 
\end{rem}

The following definition is due to Z.-J. Ruan \cite[Theorem 1.1]{Ruan1}. 
\begin{defn}
Let $(M,\Delta)$ be a locally compact quantum  group. 
We say that it is strongly Voiculescu amenable 
if there exists a net of unit vectors 
$\{\xi_j\}_{j\in \mathcal{J}}$ in $H$ such that 
for any vector $\eta$ in $H$, 
$\|W^*(\eta\otimes\xi_j)-\eta\otimes\xi_j\|$ converges to $0$.
\end{defn}

\begin{lem}\label{lem 4}
Let $(M,\Delta)$ be a locally compact quantum group.
Then the following conditions are equivalent.  
\begin{enumerate}

\item 
It is strongly Voiculescu amenable. 

\item 
There exists a net of unit vectors 
$\{\xi_j\}_{j\in \mathcal{J}}$ in $H$ with 
$\lim_j \|\lambda(\omega)\xi_{j}-\omega(1)\xi_{j}\|= 0$ 
for any functional $\omega\in M_*$. 

\item 
There exists a net of normal states $\{\omega_j\}_{j\in \mathcal{J}}$ 
on $\hat{M}$ such that 
$\{(\iota\otimes\omega_{j})(W)\}_{j\in\mathcal{J}}$ 
is a $\sigma$-weakly approximate unit of $A$.

\item 
There exists a net of normal states $\{\omega_j\}_{j\in \mathcal{J}}$ 
on $\hat{M}$ such that 
$\id : \hat{A}\longrightarrow\hat{A}$ is pointwise-weakly 
approximated by the net of unital completely positive maps 
$\{(\id\otimes\omega_j)\circ\hat{\Delta}\}_{j\in \mathcal{J}}$ 
and 
$\{(\omega_j\otimes\id)\circ\hat{\Delta}\}_{j\in \mathcal{J}}$.
  
\item 
There exists a net of normal states $\{\omega_j\}_{j\in \mathcal{J}}$ 
on $\hat{M}$ such that $\id : \hat{A}\longrightarrow\hat{A}$ 
is pointwise-norm approximated by 
the net of unital completely positive maps 
$\{(\id\otimes\omega_j)\circ\hat{\Delta}\}_{j\in \mathcal{J}}$ 
and 
$\{(\omega_j\otimes\id)\circ\hat{\Delta}\}_{j\in \mathcal{J}}$.

  \end{enumerate}

\end{lem}
\begin{proof}
 
$(1)\Rightarrow (2)$. 
It suffices to prove the statement in the case that $\omega$ 
is a normal state on $M$ by considering linear combinations. 
Since $M$ is standardly represented, 
$\omega$ is written as $\omega=\omega_\eta$ 
with a unit vector $\eta\in H$. 
For any vector $\zeta\in H$, we have 
\[
|\langle\lambda(\omega_\eta)\xi_j-\omega_\eta(1)\xi_j|\zeta\rangle|
\leq 
\|\zeta\| \|\eta\| \|W(\eta\otimes\xi_j)-\eta\otimes\xi_j\|,
\]
so the inequality 
$\|\lambda(\omega_\eta)\xi_j-\omega_\eta(1)\xi_j\| 
\leq 
\|\eta\| \|W(\eta\otimes\xi_j)-\eta\otimes\xi_j\|$ holds. 
Hence $\|\lambda(\omega_\eta)\xi_j-\omega_\eta(1)\xi_j\|$ 
converges to 0.  

$(2)\Rightarrow (3)$. 
Put $\omega_j=\omega_{\xi_j}$ for any $j$ in $\mathcal{J}$. 
Then for any operator $a\in A$ and normal functional $\theta\in M_*$, 
we have 
\begin{align*}
|\theta(a(\iota\otimes\omega_j)(W)-a)|
=&\, 
|\langle\lambda(\theta a)\xi_j-\theta(a)\xi_j|\xi_j\rangle|\\
\leq&\, 
\|\lambda(\theta a)\xi_j-(\theta a)(1)\xi_j\|.
\end{align*} 
Therefore, $a(\iota\otimes\omega_j)(W)-a$ converges to $0$ 
$\sigma$-weakly. 
Similarly we see $(\iota\otimes\omega_j)(W)a-a$ converges to $0$ 
$\sigma$-weakly.

$(3)\Rightarrow (4)$. 
Take a net of normal states $\{\omega_j\}_{j\in \mathcal{J}}$ 
of $\hat{M}_*$ which satisfies the third condition.  
Let $\omega$ be a normal functional on $M$. 
By applying Cohen's factorization theorem 
(\cite[Theorem 10, p.~61]{Bonsall and Duncan}) to 
the left $A$-module $M_\ast$, 
we get $a$ in $A$ and $\omega'$ in $M_*$ such that 
$\omega=a\omega'$. 
Then for any functional $\theta\in\hat{A}^*$, 
we have 
\begin{align*}
\theta((\omega_j\otimes\id)\circ\hat{\Delta}(\lambda(\omega)))
=&\,
\omega((\iota\otimes\theta)(W)(\iota\otimes\omega_j)(W))\\
=&\,
(\omega'(\iota\otimes\theta)(W))((\iota\otimes\omega_j)(W))a),
\end{align*} 
which converges to $\theta(\lambda(\omega))$. 
Since the linear subspace $\{\lambda(\omega);\omega\in M_*\}$ 
is norm dense in $\hat{A}$ and 
$\{\theta\circ(\omega_j\otimes\id)\circ 
\hat{\Delta}\}_{j\in \mathcal{J}}$ 
is a norm bounded family, 
$\theta((\omega_j\otimes\id)\circ\hat{\Delta}(x))$ 
converges to $\theta(x)$ for any operator $x\in\hat{A}$. 
Similarly we can see that 
$\theta((\id\otimes\omega_j)\circ\hat{\Delta}(x))$ 
converges to $\theta(x)$ for $x\in\hat{A}$.  
  
$(4)\Rightarrow (5)$. 
Take a net of normal states $\{\omega_j\}_{j\in \mathcal{J}}$ 
on $\hat{M}$ 
which satisfies the fourth condition. 
Let $\mathcal{F}$ be the set of finite subsets of $\hat{A}$. 
Take $F=\{a_1,a_2,\dots,a_k\}$ in $\mathcal{F}$ and $n$ in $\N$. 
Consider the product Banach space 
$\hat{A}_F=l_\infty$-$\sum_{x\in\mathcal{F}}\hat{A}\times\hat{A}$ 
and its dual Banach space 
$\hat{A}_F^*=l_1$-$\sum_{x\in\mathcal{F}}(\hat{A}\times\hat{A})^*$. 
Denote the following element of $\hat{A}_F$ by $x_F(\omega)$, 
\begin{align*}
\big((\omega\otimes\id)\circ\hat{\Delta}(a_1)-a_1 
&, (\id\otimes\omega)\circ\hat{\Delta}(a_1)-a_1,\\ 
(\omega\otimes\id)\circ\hat{\Delta}(a_2)-a_2 
&, 
(\id\otimes\omega)\circ\hat{\Delta}(a_2)-a_2,\\ 
&\vdots\\
(\omega\otimes\id)\circ\hat{\Delta}(a_k)-a_k 
&, 
(\id\otimes\omega)\circ\hat{\Delta}(a_k)-a_k\big).
\end{align*} 
Then $x_F(\omega_j)$ converges to $0$ weakly. 
Hence the norm closure of the convex hull of 
$\{x_F(\omega_j);j\in \mathcal{J}\}$ contains $0$. 
So there exists a normal state $\omega_{(F,n)}$ on $\hat{M}$ 
such that 
$\big\|(\omega_{(F,n)}\otimes\id)
\circ\hat{\Delta}(a)-a\big\|<{1\over n}$, 
and 
$\big\|(\id\otimes\omega_{(F,n)})
\circ\hat{\Delta}(a)-a\big\|<{1\over n}$ 
for any element $a\in F$. 
This new net $\{\omega_{(F,n)}\}_{(F,n)\in\mathcal{F}\times\N}$ 
is a desired one.

$(5)\Rightarrow (4)$. It is trivial.

$(4)\Rightarrow (3)$. Easy to prove by reversing 
the proof of $(3)\Rightarrow(4)$.

$(3)\Rightarrow (1)$. 
Take such a net of normal states $\{\omega_j\}_{j\in \mathcal{J}}$ 
on $M$. 
Since $\hat{M}$ is standardly represented, there exists a net of 
unit vectors $\{\xi_j\}_{j\in \mathcal{J}}$ in $H$ with 
$\omega_j=\omega_{\xi_j}$. 
Take a vector $\eta$ in $H$. Now by Cohen's factorization theorem, 
there exist an operator $a\in A$ and a vector $\zeta\in H$ with $\eta=a\zeta$. 
Then we have 
\begin{align*}
\|W(\eta\otimes\xi_j)-\eta\otimes\xi_j\|^2
=&\,
2\|a\zeta\|^2-2\re(\langle W(a\zeta\otimes\xi_j)|a
\zeta\otimes\xi_j\rangle)\\  
=&\,
2\|a\zeta\|^2-2\re(\langle (\iota\otimes\omega_j)(W)a
\zeta|a\zeta\rangle). 
\end{align*}
This converges to $0$. 

\end{proof}

\begin{lem}\label{lem 5}
Let $(M,\Delta)$ be a locally compact quantum  group. 
The following conditions are equivalent.  
\begin{enumerate}
 
\item 
There exists a net of unit vectors $\{\xi_j\}_{j\in \mathcal{J}}$ in $H$ 
with 
$\lim_j\|(\pi\otimes\iota)(W)^*(\eta\otimes\xi_j)-\eta\otimes\xi_j\|= 0$ 
for any representation $\{\pi,H_{\pi}\}$ of $A$ and $\eta\in H_\pi$. 

\item 
There exists a net of unit vectors $\{\xi_j\}_{j\in\mathcal{J}}$ in $H$ 
with 
$\lim_j \|(\pi\otimes\iota)(W)^\ast(\eta\otimes\xi_j)
-(\eta\otimes\xi_j)\|=0$ 
for any cyclic representation $\{\pi,H_\pi\}$ of $A$ and $\eta\in H_\pi$. 

\item 
There exists a net of unit vectors $\{\xi_j\}_{j\in \mathcal{J}}$ in $H$ 
with 
$\lim_j \|\lambda(\omega)\xi_{j}-\omega(1)\xi_{j}\|= 0$
for any functional $\omega\in A^*$. 
  
\item 
There exists a net of normal states $\{\omega_j\}_{j\in \mathcal{J}}$ 
in $\hat{M}_*$ 
such that $\{(\iota\otimes\omega_{j})(W)\}_{j\in\mathcal{J}}$ is 
a weakly approximate unit of $A$.
  
\item 
There exists a net of normal states $\{\omega_j\}_{j\in \mathcal{J}}$ 
in $\hat{M}_*$ 
such that $\{(\iota\otimes\omega_{j})(W)\}_{j\in\mathcal{J}}$ is 
a norm approximate unit of $A$.
  
\item 
There exists a net of normal states 
$\{\omega_j\}_{j\in \mathcal{J}}$ in $\hat{M}_*$ 
such that $\id:\hat{A}\longrightarrow\hat{A}$ is approximated 
in the pointwise norm topology 
by the net of unital completely positive maps 
$\{(\id\otimes\omega_j)\circ\hat{\Delta}\}_{j\in \mathcal{J}}$ and 
$\{(\omega_j\otimes\id)\circ\hat{\Delta}\}_{j\in \mathcal{J}}$.

\item 
There exists a character $\varrho$ on $\hat{A}$ with 
$(\iota\otimes\varrho)(W)=1$.
  
\item The $C^\ast$-algebra $\hat{A}$ has a character. 

\item 
There exists a state $\varrho$ on $\hat{M}$ such that 
$\varrho$ is an $\hat{A}$-linear map and 
satisfies $(\iota\otimes\varrho)(W)=1$.

\end{enumerate}
\end{lem}

\begin{proof}
$(1)\Rightarrow (2)$. It is trivial. 

$(2)\Rightarrow (3)$. 
Take such a net of unit vectors $\{\xi_j\}_{j\in \mathcal{J}}$ in $H$. 
Let $\omega$ be a state on $A$ and 
$\{H_\omega, \pi_\omega, \xi_\omega\}$ be its GNS representation. 
Then for any $\zeta$ in $H$, we have 
\[
|\langle\lambda(\omega)\xi_j-\omega(1)\xi_j|\zeta\rangle|
\leq 
\|\zeta\|\|(\pi_\omega\otimes\iota)(W)(\xi_{\omega}\otimes\xi_j)
-\xi_{\omega}\otimes\xi_j\|.
\] 
So we get 
$\|\lambda(\omega)\xi_j-\omega(1)\xi_j\|
\leq 
\|(\pi_\omega\otimes\iota)(W)(\xi_{\omega}\otimes\xi_j)
-\xi_{\omega}\otimes\xi_j\|$. 
Hence for any $\omega\in A^*$, 
$\|\lambda(\omega)\xi_j-\omega(1)\xi_j\|$ converges to $0$.

$(3)\Rightarrow (4)$. 
Put $\omega_j=\omega_{\xi_j}$ for any $j$ in $\mathcal{J}$. 
Then for any operator $a\in A$ and any functional $\theta\in A^*$, 
\[
|\theta(a(\iota\otimes\omega_j)(W)-a)|
=|\langle\lambda(\theta a)\xi_j-\theta(a)\xi_j|\xi_j\rangle| 
\leq \|\lambda(\theta a)\xi_j-(\theta a)(1)\xi_j\|.
\] 
Therefore, $a(\iota\otimes\omega_j)(W)-a$ converges to $0$ weakly. 
Similarly, we can see $(\iota\otimes\omega_j)(W)a-a$ converges to $0$ 
weakly.

$(4)\Rightarrow (5)$. 
Take such a net of normal states $\{\omega_j\}_{j\in \mathcal{J}}$ 
of $\hat{M}_*$. 
Let $\mathcal{F}$ be the set of finite subsets of $A$.
Take $F=\{a_1,a_1,\dots,a_k\}$ in $\mathcal{F}$ and $n$ in $\N$. 
Consider the product Banach space $A_F=l_\infty$-$\sum_{x\in\mathcal{F}}A$ 
and its dual Banach space $A_F^*=l_1$-$\sum_{x\in\mathcal{F}}A^*$. 
Denote the following element in $A_F$ by $x_F(\omega)$, 
\[
((\iota\otimes\omega)(W)a_1-a_1,
(\iota\otimes\omega)(W)a_2-a_2,
\dots,(\iota\otimes\omega)(W)a_k-a_k).
\] 
Now the net of the elements in $A_F$, $x_F(\omega_j)$ converges to $0$ weakly. 
Hence the norm closure of the convex hull of 
$\{x_F(\omega_j);j\in \mathcal{J}\}$ contains 0. 
So there exists a normal state of $\hat{M}_*$, $\omega_{(F,n)}$ 
such that 
$\|(\iota\otimes\omega_{(F,n)})(W)a-a\|<{1\over n}$ for any $a$ in $F$. 
Then we get a new net of normal states of $\hat{M}_*$, 
$\{\omega_{(F,n)}\}_{(F,n)\in\mathcal{F}\times\N}$ and it is 
easy to see this net is a desired one.
 
$(5)\Rightarrow (6)$. 
This is shown in a similar way to the proof of 
$(3)\Rightarrow(4)\Rightarrow(5)$ in Lemma \ref{lem 4}.
 
 $(6)\Rightarrow (7)$. 
Take such a net of normal states $\{\omega_j\}_{j\in \mathcal{J}}$ 
of $\hat{M}_*$. 
Let $\varrho$ be a weak$*$-accumulating state in 
$\hat{A}^\ast$ of the net. 
Then for any normal functional $\omega$ on $\hat{M}$, 
we have 
\begin{align*}
(\omega_j\otimes\id)(\hat{\Delta}(\lambda(\omega)))
=&\,
(\omega_j\otimes\iota\otimes\omega)(\hat{W}_{23}^\ast\hat{W}_{13}^\ast)\\ 
=&\,(\iota\otimes\omega)(\hat{W}^\ast(\omega_j\otimes\iota)(\hat{W}^\ast)).
\end{align*} 
This converges to 
$\lambda(\omega)=(\iota\otimes\omega)
(\hat{W}^\ast(\varrho\otimes\iota)(\hat{W}^\ast))$. 
Therefore we obtain $(\iota\otimes\varrho)(W)=1$ 
and easily see $\varrho$ is a character of $\hat{A}$. 

$(7)\Rightarrow (8)$. It is trivial. 

$(8)\Rightarrow (7)$. 
Take a character $\varrho$ on $\hat{A}$. Let $u$ be a unitary 
$(\iota\otimes\varrho)(W)$ in $M$. 
Then we have 
\begin{align*}
\Delta(u)
=&\,
(\iota\otimes\iota\otimes\varrho)((\Delta\otimes\iota)(W))\\
=&\,
(\iota\otimes\iota\otimes\varrho)(W_{13}W_{23})\\
=&\,u\otimes u.
\end{align*}
Therefore, for any normal functional $\omega\in M_\ast$, 
we get 
\begin{align*}
\lambda(u^\ast\omega)
=&\,
(\omega\otimes\iota)(W(u^\ast\otimes1))\\                                       =&\,
(\omega\otimes\iota)((1\otimes u^\ast)W(1\otimes u))\\
=&\, u^\ast\lambda(\omega)u.
\end{align*}
Then we obtain 
\begin{align*}
|\omega(1)|
=&\,
|(u^\ast\omega)(u)|\\ 
=&\,
|\varrho(\lambda(u^\ast\omega))|\\ 
\leq&\,
\|\lambda(u^\ast\omega)\|\\ 
=&\,
\|u^\ast\lambda(\omega)u\|\\ 
=&\,
\|\lambda(\omega)\|.
\end{align*}
Hence we can define the character $\chi$ on $\hat{A}$ with 
$\chi(\lambda(\omega))=\omega(1)$ for $\omega\in M_\ast$.

$(7)\Rightarrow (9)$. 
Take such a character $\varrho$ on $\hat{A}$ and extend it to 
the state on $\hat{M}$. Denote the extended state by $\bar{\varrho}$. 
The $\hat{A}$-linearity of $\bar{\varrho}$ easily follows from 
Stinespring's theorem.

$(9)\Rightarrow (3)$. 
Take such a state $\varrho$ on $\hat{M}$ and let 
$\{\omega_j\}_{j\in\mathcal{J}}$ be a net of normal states on 
$\hat{M}$ which weakly$*$ converges to $\varrho$. 
Then for any functional $\theta$ on $A$ and for any normal functional 
$\omega$ on $\hat{M}$, we have 
$\theta((\iota\otimes\omega_j)(W))=\omega_j((\theta\otimes\iota)(W))$. 
This converges to 
$\varrho((\theta\otimes\iota)(W))=\theta(1)$.       
  
$(3)\Rightarrow (2)$. 
Take a state $\omega$ of $A$ and let $\{H_\omega,\pi_\omega,\xi_\omega\}$ 
be the GNS representation of $\omega$. Take a vector $\eta$ in $H_\omega$. 
By Cohen's factorization theorem, there exists an operator $a\in A$ and 
a vector $\zeta\in H$ with $\eta=\pi_\omega(a)\zeta$. Then, we have 
\begin{align*} 
\|(\pi_\omega\otimes\iota)(W)\eta\otimes\xi_j-\eta\otimes\xi_j\|^2
=&\,
2\|\pi_\omega(a)\zeta\|^2\\
&
-\re(\langle(\pi_\omega\otimes\iota)(W)\pi_\omega(a)\zeta\otimes\xi_j|
\pi_\omega(a)\zeta\otimes\xi_j\rangle)\\ 
=&\,
2\|\pi_\omega(a)\zeta\|^2-\re(\langle
\pi_\omega(a^\ast(\iota\otimes\omega_j)(W)a)\zeta|\zeta\rangle). 
\end{align*}
This converges to $0$. 

$(2)\Rightarrow (1)$. 
Let $\{\pi,H_\pi\}$ be a representation of $A$. 
We may assume that this representation is nondegenerate. 
Then $\pi$ is decomposed to the direct sum of cyclic representation. 
This observation derives the statement of 1. 
                   
\end{proof}

By the previous two lemmas, 
we obtain the following result. 

\begin{cor}\label{cor 3}
Let $(M,\Delta)$ be a locally compact quantum  group. 
Then the following statements are equivalent.

\begin{enumerate}

\item 
It is strongly Voiculescu amenable. 

\item 
There exists a net of unit vectors $\{\xi_j\}_{j\in \mathcal{J}}$ in $H$ 
satisfying 
$\lim_j\|(\pi\otimes\iota)(W)^*(\eta\otimes\xi_j)-\eta\otimes\xi_j\|=0$
for any representation $\{\pi,H_{\pi}\}$ of $A$ and $\eta\in H_\pi$.

\end{enumerate}
\begin{proof}

$(2)\Rightarrow (1)$. It is trivial.

$(1)\Rightarrow (2)$. 
The latter statement is equivalent to 
the statement (6) in Lemma \ref{lem 5} 
and it is the same as 
the statement (5) 
in Lemma \ref{lem 4} which is equivalent to strong Voiculescu amenability. 
\end{proof}
\end{cor}

The following results have been already known in various settings 
\cite[Theorem 4.2]{Bedos Murphy Tuset1}.

\begin{cor}
Let $(M,\Delta)$ be a locally compact quantum  group. 
Then the following statements are equivalent.

\begin{enumerate}

\item It is strongly Voiculescu amenable.

\item 
The $C^*$-algebra $\hat{A}$ has a character $\varrho$ 
with $(\iota\otimes\varrho)(W)=1$.

\item The $C^*$-algebra $\hat{A}$ has a character. 

\end{enumerate}
\end{cor}

We begin to treat the discrete quantum groups from now. 
The discreteness of $(M,\De)$ is always assumed. 
The following theorem is the main result of this paper.

\begin{thm}\label{main}
If $(M,\Delta)$ is a discrete quantum group, 
the following statements are equivalent. 
\begin{enumerate}
\item It has an invariant mean.
\item It is strongly Voiculescu amenable. 
\item The $C^*$-algebra $\hat{A}$ is nuclear and has a character.
\item 
The von Neumann algebra $\hat{M}$ is injective and has an 
$\hat{A}$-linear state. 
\end{enumerate}
\end{thm}

We remark that the equivalence $(2)\Leftrightarrow(3)\Leftrightarrow(4)$ 
have been already proved by the combination of 
\cite[Theorem 5.2]{Bedos Conti Tuset}, 
\cite[Corollary 2.9]{Bedos Murphy Tuset1}, 
\cite[Theorem 4.8]{Bedos Murphy Tuset2}, 
\cite[Theorem 1.1, Theorem 4.2, Corollary 4.3]{Bedos Murphy Tuset3}, 
however we give a different proof of deriving nuclearity by using 
completely positive approximation property. 
Before proving it, we state the following corollary, 
where a subgroup means a $C^*$-algebraic 
compact quantum group $(C(K),\Delta_K)$ with a compact group $H$ which 
is an image space of $*$-homomorphism $r:\hat{A}\longrightarrow C(K)$ with 
$\Delta_K \circ r=(r\otimes r)\circ \hat{\Delta}$. 

\begin{cor}
Let $(\hat{A},\hat{\Delta})$ be a compact quantum group 
which has a subgroup. 
Then the dual discrete quantum group $(A,\Delta)$ is amenable.  
\end{cor}

A character is constructed by the composition of 
the restriction map and a character on continuous function algebra of 
the subgroup. 
Therefore, 
if a compact quantum group is deformed by some parameters 
from a ordinary group 
and a subgroup (e.g. maximal torus) becomes a non-deformed 
subgroup, 
then its dual discrete quantum group is amenable. 
Let $\M=\oplus_{\al\in I}Mz_\al$ 
be a matrix decomposition as in Section 2, 
where $Mz_\alpha\cong \Ma_{n_\al}(\C)$. 
For a finite subset $F$ in $I$, 
let us denote  $z_F=\sum_{\alpha\in F}z_\alpha$. 
Note that $z_F H=\Lambda(Mz_F)$ is finite dimensional subspace. 

\begin{lem}\label{lem 6}
If $F$ is a finite subset of $I$, then the linear subspace 
$K_F=\{\lambda(\omega z_F);\\\omega\in M_*\}\subset\hat{A}$ is 
finite dimensional. 
\end{lem}

\begin{proof}
Take a normal functional $\omega$ in $\mathcal{I}$ as introduced in Section 2 
and an operator $x$ in $n_\varphi$. Then we have 
\begin{align*}
  \omega(z_F x^*)=&\,\omega((z_F x)^*)\\
  =&\,\langle\lambda(\omega)|z_F\Lambda(x)\rangle\\
  =&\,\langle z_F\lambda(\omega)\xi_{\hat{\varphi}}|\Lambda(x)\rangle.
\end{align*}   
So $\omega z_F$ is in $\mathcal{I}$ and we obtain 
$\lambda(\omega z_{F})\xi_{\hat{\varphi}}
=z_F\lambda(\omega)\xi_{\hat{\varphi}}\in z_FH$. Since $\xi_{\hat{\varphi}}$ is 
a separating vector for $\hat{M}$ and $\mathcal{I}$ is norm dense in $M_*$, 
the statement follows. 
\end{proof}

We give a proof of Theorem \ref{main} partly first. 

{\it Proof of Theorem \ref{main}} 
($(2)\Rightarrow(3)\Rightarrow(4)\Rightarrow(1)$). 
Fix a complete orthonormal system 
$\{e_p\}_{p\in\mathcal{P}}$ of $H$. 

$(2)\Rightarrow(3)$. 
We show that $\hat{A}$ has completely positive approximation property. 
Set $T_\omega=(\iota\otimes\omega)\circ\hat{\Delta}$ for 
$\omega\in\hat{M}_*$. By assumption, there exists 
a net of normal states $\{\omega_j\}_{j\in \mathcal{J}}$ 
on $\hat{M}$ 
satisfying the statement (6) in Lemma \ref{lem 5}, 
that is, $T_{\omega_j}$ converges to the identity map of $\hat{A}$ 
in the pointwise norm topology. 
Since $\hat{M}$ is standardly represented, 
each $\omega_j$ is a vector 
state defined by a unit vector $\xi_j\in H$. 
Take a finite subset $F$ of $I$.
Then for any operator $x\in\hat{A}$, we have 
\begin{align*}
T_{\omega_{z_F\xi_j}}(x)
=&\,
(\iota\otimes\omega_{z_F\xi_j})(\hat{\Delta}(x))\\ 
=&\,
\sum_{p\in\mathcal{P}}{(\iota\otimes\omega_{z_F\xi_{j},\,e_p})(\hat{W})}^*
(\iota\otimes\omega_{z_F\xi_{j},\,e_p})((1\otimes x)\hat{W})\\
=&\,
\sum_{p\in\mathcal{P}}\lambda((\omega_{e_p,\,\xi_{j}})z_F)
\lambda((\omega_{x^* e_p,\,\xi_{j}})z_F)^*.
\end{align*}
Since $\lambda((\omega_{e_p,\,\xi_{j}})z_F)
\lambda((\omega_{x^\ast e_p,\,\xi_{j}})z_F)^*$ is in the finite dimensional 
linear subspace $K_F K_F^\ast=\spa\{ab^\ast;a,b\in K_F\}\subset\hat{A}$, 
so is $T_{\omega_{z_F\xi_j}}(x)$. 
Hence the completely positive map $T_{\omega_{z_F\xi_j}}$ has finite rank. 
For $n\in \N$ and $j\in\mathcal{J}$, 
take an finite subset $F(n,j)$ of $I$ with 
$\|T_{\omega_{z_{F(n,j)}\xi_j}}-T_{\omega_{\xi_j}}\|<{1\over n}$. 
Then we get a new net of completely positive maps  
$\{T_{\omega_{z_{F(n,j)}\xi_j}}\}_{(n,j)\in\N \times\mathcal{J}}$
of finite rank. 
This net gives a completely positive approximation of the identity map of 
$\hat{A}$, 
hence $\hat{A}$ is a nuclear $C^*$-algebra.
                  
$(3)\Rightarrow(4)$. 
nuclearity $\hat{A}$ implies injectivity $\hat{M}$. 
Extend the character on $\hat{A}$ to the state on $\hat{M}$. 
We can easily see the $\hat{A}$-linearity of this state 
by Stinespring's theorem. 

$(4)\Rightarrow(1)$. 
Let $\varrho$ be an $\hat{A}$-linear state on $\hat{M}$. 
There exists a conditional expectation $E$ from $B(H)$ onto $\hat{M}$. 
We set a state $m$ on $M$ by $m=\varrho\circ E|_M$. 
Then for any vector $\xi\in H$ and for any operator $x\in M$, we have 
\begin{align*}
\omega_{\xi}*m(x)
=&\,
m((\omega_{\xi}\otimes\iota)(\Delta(x)))\\
=&\,
m\Big(\sum_{p\in\mathcal{P}}
(\omega_{\xi,\,e_p}\otimes\iota)(W)^* x
(\omega_{\xi,\,e_p}\otimes\iota)(W)\Big)\\
=&\,
\sum_{p\in\mathcal{P}}m\big(
(\omega_{\xi,\,e_p}\otimes\iota)(W)^* x
(\omega_{\xi,\,e_p}\otimes\iota)(W)\big)\\
=&\,
\sum_{p\in\mathcal{P}}\varrho\big(\big(
(\omega_{\xi,\,e_p}\otimes\iota)(W)^* E(x)
(\omega_{\xi,\,e_p}\otimes\iota)(W)\big)\\
=&\,
\sum_{p\in\mathcal{P}}
\varrho((\omega_{\xi,\,e_p}\otimes\iota)(W)^*) 
\varrho(E(x))
\varrho((\omega_{\xi,\,e_p}\otimes\iota)(W))\\ 
=&\,
\sum_{p\in\mathcal{P}}
\varrho((\omega_{\xi,\,e_p}\otimes\iota)(W)^*) 
\varrho((\omega_{\xi,\,e_p}\otimes\iota)(W))
\varrho(E(x))\\
=&\,
\sum_{p\in\mathcal{P}}
\varrho((\omega_{\xi,\,e_p}\otimes\iota)(W)^*
(\omega_{\xi,\,e_p}\otimes\iota)(W))
\varrho(E(x))\\ 
=&\,
\varrho(\omega_{\xi}(1))\varrho(E(x))\\
=&\,
\omega_{\xi}(1)m(x),
\end{align*}
where we have used the norm convergence of 
$\sum_{p\in\mathcal{P}}
(\omega_{\xi,\,e_p}\otimes\iota)(W)^* x
(\omega_{\xi,\,e_p}\otimes\iota)(W)$ in the third equality. 
Therefore, $m$ is a left invariant mean on $M$.\hfill$\square$

Now we are going to prove the implication $(1)\Rightarrow(2)$ 
of Theorem \ref{main} after proving the several lemmas. 
As usual, 
$L^\infty(\R)$ means the von Neumann algebra 
which consists of essentially bounded measurable functions 
with respect to the Lebesgue measure. 

\begin{lem}\label{lem 7}
Let $m_{\R}$ be an invariant mean of $L^\infty(\R)$. For any $\omega$ in $M_*$ define 
the $\sigma^\vph$-invariant functional $\omega'$ by 
$\omega'(x)=m_{\R}(\{t\mapsto\omega(\si_t^\vph(x))\})$ 
for x in $M$. Then $\omega'$ is a normal functional with 
$\|\omega'\|\leq\|\omega\|$.    
\end{lem}
\begin{proof}
For $\omega$ in $M_*$ set 
$f_{\omega,\,x}(t)=\{t\mapsto\omega(\si_t^\vph(x))\}$ in $C^b(\R)$. 
For any finite subset $F$ in $I$, 
$|f_{\omega,\,x}(t)-f_{\omega z_F,\,x}(t)| 
\leq \|\omega-\omega z_F\| \|x\|$, 
hence 
$\|\omega'-(\omega z_F)'\| \leq \|\omega -\omega z_F\|$. 
By the normality of $\omega$, 
$\lim_F (\omega z_F)'=\omega'$. Notice that $Mz_F$ is finite dimensional, so 
$(\omega z_F)'=(\omega z_F)'z_F$ is a normal functional. 
\end{proof}

Note that this averaging procedure can be also done by considering 
a conditional expectation from $M$ to $M_\vph$. 
Recall that 
we have fixed a matrix unit 
$\{e(\alpha)_{kl}\}_{1\leq k,\,l\leq n_\alpha}$ of 
$Mz_\alpha$ for each $\alpha$ in $I$, such that they are diagonalizing 
$h z_\alpha$ as 
$h z_\alpha
=\sum_{k=1}^{n_\alpha} \nu(\alpha)_k e(\alpha)_{kk}$, 
where $\nu(\alpha)_k$ denotes a positive real number. 

\begin{lem}\label{lem 8}
If $\M$ has an invariant mean $m$, there exists a net of normal states 
$\{\omega_j\}_{j\in \mathcal{J}}$ on $M$, which satisfies the following 
two conditions.
\begin{enumerate}

\item 
$\lim_j\|\omega*\omega_j -\omega(1)\omega_j\|=0$ 
for any $\omega$ in $M_*$. 
\item 
$\omega_j\circ\si_t^\vph=\omega_j$ for any $t\in \R$. 

\end{enumerate}
\end{lem}

\begin{proof}
Firstly we show the existence of a net satisfying the first condition. 
Since the convex hull of the vector states is weak$*$-dense in the state 
space of $M$, there exists a net of normal states $\{\chi_j\}_{j\in \mathcal{J}}$ in 
$M_*$ such that $m=\mbox{w*-}\lim_j \chi_j$. 
Let $\mathcal{F}$ be the set of finite subsets of $M_*$.
For $F=\{\omega_1,\omega_2,\dots,\omega_k\}\in\mathcal{F}$, consider the Banach space 
${(M_*)}_F=l_1$-$\sum_{\omega\in F}{M_*}$ and its dual Banach space 
$M_F=l_\infty$-$\sum_{\omega\in F}M$. 
Set $$x_F(\chi)=(\omega_1*\chi-\omega_1(1)\chi, 
                       \omega_k*\chi-\omega_k(1)\chi,\dots,
                       \omega_k*\chi-\omega_k(1)\chi),$$ for $\chi$ in $M_*$. 
Then $x_F(\chi_j)$ converges to 0 weakly. 
So the norm closure of the convex hull of $\{x_F(\chi_j);j\in \mathcal{J}\}$ contains $0$. 
Hence for any $n$ in $\N$, there exists $\chi_{(F,\,n)}$ such that 
$\|\omega*\chi_{(F,\,n)}-\omega(1)\chi_{(F,\,n)}\|<{1\over n}$ for an $\omega$ in $F$. 
The new net $\{\chi_{(F,\,n)}\}_{(F,\,n)\in\mathcal{F}\times\N}$ 
is a desired one. 
Next we show existence of a net satisfying the both conditions. 
Let $\{\omega_j\}_{j\in \mathcal{J}}$ be a net satisfying the first condition. 
By the previous Lemma, $\omega_j'$ is normal. 
We show that the net $\{\omega'_j\}_{j\in \mathcal{J}}$ 
satisfies the first condition. 
For 
$\omega=\omega_{\Lambda(e(\alpha)_{kl}),\,\Lambda(e(\alpha)_{mn})}$, 
we have 
$\omega\circ\si_{-t}^\vph
=\nu_{k}(\alpha)^{it}\nu_{l}(\alpha)^{-it}
\nu_{m}(\alpha)^{-it}\nu_{n}(\alpha)^{it}\omega$. 
Then we obtain 
\begin{align*}
&\,|\omega*\omega'_j(x)-\omega(1)\omega'_j(x)|\\
=&\,
|m_\R(\{t\mapsto (\omega\otimes\omega_j\circ\si_t^\vph)(\Delta(x))
-\omega(1)\omega_j(\si_t^\vph(x))\})|\\
=&\,
|m_\R(\{t\mapsto (\omega\circ\si_{-t}^\vph \otimes\omega_j)
(\Delta(\si_t^\vph(x)))
-\omega(\si_{-t}^\vph(1))\omega_j(\si_t^\vph(x))\})|\\
\leq&\, 
\sup_{t\in\R}\{|(\omega\circ\si_{-t}^\vph)*\omega_j(\si_t^\vph(x))
-\omega(\si_{-t}^\vph(1))\omega_j(\si_t(x))|\} \\
  \leq&\, \sup_{t\in\R}\{\|(\omega\circ\si_{-t}^\vph)*\omega_j
  -\omega\circ \si_{-t}^\vph(1)\,\omega_j\|\} \|x\|\\ 
  =&\,\sup_{t\in\R}\{\|\nu_{k}(\alpha)^{it}\nu_{l}(\alpha)^{-it}
                           \nu_{m}(\alpha)^{-it}\nu_{n}(\alpha)^{it}
                          (\omega*\omega_j-\omega(1)\omega_j)\|\}\|x\|\\
  \leq&\, \|\omega*\omega_j-\omega(1)\omega_j\| \|x\|. 
  \end{align*}
So for this $\omega$, we have 
$\|\omega*\omega'_j-\omega(1)\omega'_j\| \leq \|\omega*\omega_j-\omega(1)\omega_j\|$ 
and therefore $\|\omega*\omega'_j-\omega(1)\omega'_j\|$ converges to 0. 
By taking the linear combination for $\omega$, we see that 
$\|\omega*\omega'_j-\omega(1)\omega'_j\|$ converges to 0 for any normal functional 
with $\omega=\omega z_\alpha$. 
Take a normal functional $\omega$ and a positive $\varepsilon$. 
Then there exists a finite subset $F$ of $I$ and $j_0$ in $\mathcal{J}$ such that 
$\|\omega z_F-\omega\|<\varepsilon$ and 
$\|\omega z_F*\omega'_j-\omega z_F(1)\omega'_j\|<\varepsilon$ for $j\geq j_0$. 
Then we have 
\begin{align*}
&\,\|\omega*\omega'_j-\omega(1)\omega'_j\|\\
\leq&\, \|(\omega-\omega z_F)*\omega'_j\|
         +\|\omega z_F*\omega'_j-\omega z_F(1)\omega'_j\|
                    +\|(\omega z_F(1)-\omega(1))\omega'_j\|\\
<&\, \varepsilon+\varepsilon+\varepsilon=3\varepsilon,
\end{align*}
for $j\geq j_0$. 
Therefore, $\|\omega*\omega'_j-\omega(1)\omega'_j\|$ converges to 0.  
\end{proof}

\begin{lem}\label{lem 9}
If $\M$ has an invariant mean, there exists a net of unit vectors $\{\xi_j\}_{j\in \mathcal{J}}$ in 
H which satisfies the following four conditions. 
\begin{enumerate}

\item 
$\lim_j\|\omega*\omega_{\xi_j}-\omega(1)\omega_{\xi_j}\|=0$ for any 
$\omega$ in $M_*$. 

\item $\vDe_\vph^{it}\xi_j=\xi_j$ for any $t$ in $\R$. 

\item For any $j$ in $\mathcal{J}$, 
there exists a finite subset $F_j$ of $I$ with $z_{F_j}\xi_j=\xi_j$. 

\item 
For any $j$ in $\mathcal{J}$, the vector $\xi_j$ is in the convex cone 
$\mathcal{P}_{\varphi}^\natural$. 
\end{enumerate}
\end{lem}

\begin{proof}
There exists a net of normal state $\{\omega_j\}_{j\in \mathcal{J}}$ which 
satisfies the two conditions of the previous lemma. 
If necessary, by cutting and normalizing, 
we may assume that there exists a finite subset $F_j$ of $I$ 
 such that $\omega_j z_{F_j}=\omega_j$ for any $j$ in $\mathcal{J}$. 
The von Neumann algebra $M$ is standardly represented, so there 
exists a unique net of unit vectors $\{\xi_j\}_{j\in \mathcal{J}}$ 
in $\mathcal{P}_{\varphi}^\natural$ such that $\omega_j=\omega_{\xi_j}$ and 
$z_{F_j}\xi_j=\xi_j.$ 
From the assumption 
$\omega_j=\omega_{j}\circ\si_{-t}^\vph=\omega_{\vDe_\vph^{it}\xi_j}$, 
the uniqueness of $\xi_j$ implies $\vDe_\vph^{it}\xi_j=\xi_j$. 
\end{proof}

Let $\{\xi_j\}_{j\in \mathcal{J}}$ be a net of unit vectors in $H$ in the previous Lemma. 
Since $z_{F_j}H=\Lambda(Mz_{F_j})$, there exists $x_j$ in $M$ such that 
$\xi_j=\Lambda(x_j)$. The operator 
$x_j=x_j z_{F_j}$ is in $M_\varphi$ and satisfies 
$\varphi(x_j^* x_j)=1$ and $x_j=x_j^*$. 
We prepare some notations. 
For an operator $X$ in $M\otimes M$, $X(\alpha)$ means 
the operator $X(z_\alpha\otimes1)$ in $Mz_\alpha\otimes M$. 
An operator $Y$ in $Mz_\alpha\otimes M$ is written as 
$Y=\sum_{k,\,l=1}^{n_\alpha}e(\al)_{kl}\otimes Y_{kl}$, 
where $\{Y_{kl}\}_{kl}$ are operators in $M$. 

\begin{lem}\label{lem 10}
Let $x=x^*$ be in $M_\varphi\cap M_{z_{F}}$, where $F$ is a finite subset of $I$, 
and $\alpha$ be in $I$. Then the following inequality holds. 
\begin{align*}
&\|\omega_{\hat{J}\Lambda(e(\alpha)_{k1}),\,\hat{J}\Lambda(e(\alpha)_{l1})}
      *\omega_{\Lambda(x)}
   -\omega_{\hat{J}\Lambda(e(\alpha)_{k1}),\,\hat{J}\Lambda(e(\alpha)_{l1})}(1)
       \omega_{\Lambda(x)}\|\\ 
&\hspace{30pt}\geq \nu(\alpha)_k^{-{1\over2}}\nu(\alpha)_l^{1\over2} 
                   \varphi(e(\alpha)_{11})\varphi(|X(\alpha)_{kl}|),
\end{align*} 
where $X=\Delta(x^2)-(1\otimes x^2)$.  
\end{lem}

\begin{proof}
We simply write $e_{kl}$, $\nu_k$ and $X_{kl}$ 
for $e(\alpha)_{kl}$, $\nu(\alpha)_k$ and $X(\alpha)_{kl}$ respectively. 
Since $X(\alpha)=\sum_{1\leq k,\,l\leq n_\alpha}(e_{kl}\otimes X_{kl})$ is in 
$(M\otimes M)_{\varphi\otimes\varphi}$ 
and 
$\sigma_t^\varphi(e_{kl})=\nu_k^{it}\nu_l^{-it}e_{kl}$, 
we have 
$\sigma_t^\varphi(X_{kl})=\nu_k^{-it}\nu_l^{it}X_{kl}$. 
Let $X_{kl}=v_{kl}|X_{kl}|$ be the polar decomposition of $X_{kl}$. 
Put $a_{kl}=v_{kl}^*$. 
Then $\sigma_t^\varphi(a_{kl})=\nu_k^{it}\nu_l^{-it}a_{kl}$. Then we have  
\begin{align*}
&
(\omega_{\hat{J}\Lambda(e_{k1}),\,\hat{J}\Lambda(e_{l1})}*\omega_{\Lambda(x)}
-\omega_{\hat{J}\Lambda(e_{k1}),\,\hat{J}\Lambda(e_{l1})}(1)
\omega_{\Lambda(x)})(a_{kl})\\ 
=\,&
\langle \Delta(a_{kl})(\hat{J}\Lambda(e_{k1})\otimes\Lambda(x))|
\hat{J}\Lambda(e_{l1})\otimes\Lambda(x)\rangle
-\langle\Lambda(e_{l1})|\Lambda(e_{k1})\rangle
\langle a_{kl}\Lambda(x)|\Lambda(x)\rangle\\
=\,&
\langle(1\otimes a_{kl})
W(\hat{J}\otimes J)(\Lambda(e_{k1})\otimes\Lambda(x))|
W(\hat{J}\otimes J)(\Lambda(e_{l1})\otimes\Lambda(x))\rangle\\
&\hspace{30pt}
-\langle\Lambda(e_{l1})|\Lambda(e_{k1})\rangle
\langle a_{kl}\Lambda(x)|\Lambda(x)\rangle\\
=\,&
\langle(1\otimes Ja_{kl}^* J)W^*(\Lambda(e_{l1})\otimes\Lambda(x))|
W^*\Lambda(e_{k1})\otimes\Lambda(x)\rangle\\
&\hspace{30pt}
-\langle\Lambda(e_{l1})|\Lambda(e_{k1})\rangle
\langle Ja_{kl}^* J\Lambda(x)|\Lambda(x)\rangle\\
=\,&
\langle(\Lambda\otimes\Lambda)(\Delta(x)(e_{l1}\otimes1)
(1\otimes(\nu_k^{1\over2}\nu_l^{-{1\over2}}a_{kl})))|
(\Lambda\otimes\Lambda)(\Delta(x)(e_{k1}\otimes1))\rangle\\
&\hspace{30pt}
-\langle\Lambda(e_{l1})|\Lambda(e_{k1})\rangle
\langle\Lambda(x\nu_k^{1\over2}\nu_l^{-{1\over2}}a_{kl})|\Lambda(x)\rangle\\
=\,&
\nu_k^{1\over2}\nu_l^{-{1\over2}}\{(\varphi\otimes\varphi)((e_{1k}\otimes1)
\Delta(x^2)(e_{l1}\otimes1)(1\otimes a_{kl}))\\
&\hspace{30pt}
-(\varphi\otimes\varphi)((e_{1k}\otimes1)(1\otimes x^2)(e_{l1}\otimes1)
(1\otimes a_{kl}))\}\\
=\,&
\nu_k^{1\over2}\nu_l^{-{1\over2}}(\varphi\otimes\varphi)
(e_{11}\otimes X_{kl}a_{kl})\\
=\,&
\nu_k^{1\over2}\nu_l^{-{1\over2}}\varphi(e_{11})\varphi(X_{kl}a_{kl})\\
=\,&
\nu_k^{1\over2}\nu_l^{-{1\over2}}\nu_k^{-1}\nu_l\varphi(e_{11})
\varphi(a_{kl}X_{kl})\\
=\,&
\nu_k^{-{1\over2}}\nu_l^{1\over2}\varphi(e_{11})\varphi(|X_{kl}|). 
\end{align*}
Therefore, we obtain 
\begin{align*}
&\|\omega_{\hat{J}\Lambda(e_{k1}),\,\hat{J}
\Lambda(e_{l1})}*\omega_{\Lambda(x)}
-\omega_{\hat{J}\Lambda(e_{k1}),\,\hat{J}\Lambda(e_{l1})}(1)
\omega_{\Lambda(x)}\|\\
&\hspace{30pt}
\geq\nu_k^{-{1\over2}}\nu_l^{1\over2}\varphi(e_{11})\varphi(|X_{kl}|).
\end{align*}
\end{proof}

\begin{lem}\label{lem 11}
Let $N$ be a von Neumann algebra and $\theta$ be an n.s.f. weight on 
$N$. Let $\Ma_n(\C)$ be a matrix algebra with the matrix unit 
$\{e_{ij}\}_{1\leq i,\,j\leq n}$ and 
$\chi=\Tr_{h}$ be an n.s.f. weight on $\Ma_n(\C)$ with 
$h=\sum_{i=1}^{n}\lambda_i e_{ii}$, $\lambda_{i}>0$. 
If $\{A_j\}_{j\in \mathcal{J}}$ is a net in 
$(\Ma_n(\C)\otimes N)_{\chi\otimes\theta}$ such that 
$\theta(|(A_j)_{kl}|)<\infty$ and $\lim_{j}\theta(|(A_j)_{kl}|)=0$ for any 
$k,l=1,2,\dots,n$, 
then 
$\lim_{j}\theta(|A_j|_{kl})=0$ for any $k,l=1,2,\dots,n$.  
\end{lem}
\begin{proof}
Let $A_j=V_j|A_j|$ be the polar decomposition in 
$(\Ma_n(\C)\otimes N_\theta)_{\chi\otimes\theta}$. 
Then for each $k,l$, 
we obtain $\sigma_t^\theta((V_j)_{kl})=\lambda_k^{-it}\lambda_l^{it}(V_j)_{kl}$ and 
$\sigma_t^\theta((A_j)_{kl})=\lambda_k^{-it}\lambda_l^{it}(A_j)_{kl}$. 
Since $|A_j|_{kl}=\sum_{m=1}^{n}(V_j^*)_{km} (A_j)_{ml}$ and each $(V_j)_{km}$ is 
analytic, 
$\theta(|A_j|_{kl})$ is well-defined. 
Let $(A_j)_{ml}=v_{j,m,l}|(A_j)_{ml}|$ be 
the polar decomposition. 
Then we have
\begin{align*}
|\theta(|A_j|_{kl})|\leq&\,\sum_{m=1}^{n}|\theta((V_j^*)_{km} (A_j)_{ml})|\\
=&\,\sum_{m=1}^{n}
\big|\big\langle\Lambda_{\theta}\big(|(A_j)_{ml}|^{1\over2}\big)
      \big|\Lambda_{\theta}\big(|(A_j)_{ml}|^{1\over2}v_{j,m,l}^\ast 
           (V_j)_{mk}\big)\big\rangle\big|\\
=&\,\sum_{m=1}^{n}
\big|\big\langle\Lambda_{\theta}\big(|(A_j)_{ml}|^{1\over2}\big)
      \big|J_{\theta}\sigma_{i\over2}^\theta(v_{j,m,l}^\ast (V_j)_{mk})^\ast J_{\theta}
           \Lambda_{\theta}\big(|(A_j)_{ml}|^{1\over2}\big)\big\rangle\big|\\
\leq&\,\sum_{m=1}^{n}\|\sigma_{i\over2}^\theta(v_{j,m,l}^\ast (V_j)_{mk})\|
        \theta(|(A_j)_{ml}|)\\
=&\,\lambda_k^{-{1\over2}}\lambda_l^{1\over2}\sum_{m=1}^{n}
        \theta(|(A_j)_{ml}|).
\end{align*}
Hence $\lim_j\theta(|A_j|_{kl})=0$. 
\end{proof}

\begin{lem}\label{lem 12}
Let $x=x^*$ be in $M_\varphi\cap M_{z_{F}}$, where $F$ is a finite subset of $I$ 
such that $\varphi(x^2)=1$ and $\Lambda(x)$ is in $\mathcal{P}_\varphi^\natural$. 
Let $\alpha$ be in $I$. 
Then the following inequality holds.
\begin{align*}
&\big\|W^*\big(\Lambda(e(\alpha)_{k1})\otimes\Lambda(x)\big)
          -\Lambda(e(\alpha)_{k1})\otimes\Lambda(x)\big\|^2\\
&\hspace{20pt}\leq 2\max_{1\leq k \leq n_\alpha}\{\nu(\alpha)_k^{-1}\nu(\alpha)_1\}\cdot
\varphi(z_\alpha)^{1\over2}(\varphi\otimes\varphi)(|X(\alpha)|)^{\frac{1}{2}},
\end{align*}
where 
$X=\Delta(x^2)-1\otimes x^2$.
\end{lem}
\begin{proof}
We use the notations in Lemma \ref{lem 10}. We have 
\begin{align*}
&\big\|W^*\big(\Lambda(e_{k1})\otimes\Lambda(x)\big)
          -\Lambda(e_{k1})\otimes\Lambda(x)\big\|^2\\
=\,&2\varphi(e_{11})-2\re\big\langle W^*\big(\Lambda(e_{k1})\otimes\Lambda(x)\big)|
                              \Lambda(e_{k1})\otimes\Lambda(x)\big\rangle\\
=\,&2\varphi(e_{11})-2\re(\varphi\otimes\varphi)((e_{1k}\otimes1)
                              \Delta(x)(e_{k1}\otimes x^*))\\
=\,&2\varphi(e_{11})-2\re\nu_k^{-1}\nu_1(\varphi\otimes\varphi)((e_{kk}\otimes1)
                              (1\otimes x^*)\Delta(x))\\
=\,&2\nu_k^{-1}\nu_1\re\{\varphi(e_{kk})-(\varphi\otimes\varphi)((e_{kk}\otimes1)
                              (1\otimes x^*)\Delta(x))\}\\
=\,&2\nu_k^{-1}\nu_1\re\{(\varphi\otimes\varphi)\big((e_{kk}\otimes1)(1\otimes x^*)
                              (1\otimes x-\Delta(x))\big)\}\\
\leq\,& 2\max_{1\leq k \leq n_\alpha}\{\nu_k^{-1}\nu_1\}\cdot
      \sum_{k=1}^{n_\alpha}
      \re\{(\varphi\otimes\varphi)\big((e_{kk}\otimes x^*)
                              (1\otimes x-\Delta(x))\big)\}\\
=\,&2\max_{1\leq k \leq n_\alpha}\{\nu_k^{-1}\nu_1\}\cdot
      \re\{(\varphi\otimes\varphi)\big((z_\alpha\otimes x^*)
                              (1\otimes x-\Delta(x))\big)\}\\
\leq\,& 2\max_{1\leq k \leq n_\alpha}\{\nu_k^{-1}\nu_1\}\cdot
      |(\varphi\otimes\varphi)\big((z_\alpha\otimes1)(1\otimes x^*)
                              (1\otimes x-\Delta(x))\big)|\\
\leq\,& 2\max_{1\leq k \leq n_\alpha}\{\nu_k^{-1}\nu_1\}\cdot
    (\varphi\otimes\varphi)(z_\alpha\otimes x^* x)^{1\over2}\\  
&\hspace{30pt}\cdot(\varphi\otimes\varphi)\big((z_\alpha\otimes1)
                (1\otimes x-\Delta(x))^*(1\otimes x-\Delta(x))
\big)^{\frac{1}{2}}\\
=\,&2\max_{1\leq k \leq n_\alpha}
\{\nu_k^{-1}\nu_1\}\cdot\varphi(z_\alpha)^{\frac{1}{2}}
\big\|(\Lambda\otimes\Lambda)(z_\alpha\otimes x)
-(\Lambda\otimes\Lambda)((z_\alpha\otimes1)\Delta(x))\big\|.
\end{align*}
From the assumption: 
$\Lambda(x)=z_F\Lambda(x)\in\mathcal{P}_{\varphi}^\natural$, 
there exists a sequence $\{y_n\}_{n\in\N}\in Mz_F$ with 
$\lim_{n}y_n J\Lambda(y_n)=\Lambda(x)$. Then we obtain 
\begin{align*}
(\Lambda\otimes\Lambda)((z_\alpha\otimes1)\Delta(x))=&\,
W^\ast(\Lambda\otimes\Lambda)(z_{\alpha}\otimes x)\\
=&\,\lim_n W^\ast(\Lambda\otimes\Lambda)(z_{\alpha}\otimes 
                  y_n\sigma_{i\over2}^\varphi(y_n)^\ast)\\
=&\,\lim_n (\Lambda\otimes\Lambda)
     \big(\Delta(y_n)(\sigma_{i\over2}^\varphi\otimes\sigma_{i\over2}^\varphi)
         (\Delta(y_n))^*
       (z_\alpha\otimes1)\big)\\
=&\,\lim_n \Delta(y_n)(z_\alpha\otimes1)(J\otimes J)(\Lambda\otimes\Lambda)
         (\Delta(y_n)(z_\alpha\otimes1)).
\end{align*}
Therefore, we see $(\Lambda\otimes\Lambda)((z_\alpha\otimes1)\Delta(x))\in 
\mathcal{P}_{\varphi\otimes\varphi}^\natural$. 
By using the Powers-St\o rmer inequality 
(see, for example, \cite{Haagerup}, \cite{Powers-Stomer}), we obtain 
\begin{align*}
&\,\|(\Lambda\otimes\Lambda)(z_\alpha\otimes x)
            -(\Lambda\otimes\Lambda)((z_\alpha\otimes1)\Delta(x))\|\\
\leq&\,\|\omega_{(\Lambda\otimes\Lambda)(z_\alpha\otimes x)}
-\omega_{(\Lambda\otimes\Lambda)
((z_\alpha\otimes1)\Delta(x))}\|^{\frac{1}{2}}\\
=&\,(\varphi\otimes\varphi)(|z_\alpha\otimes x^2
           -(z_\alpha\otimes1)\Delta(x^2)|)^{\frac{1}{2}}\\
=&\,(\varphi\otimes\varphi)(|X(\alpha)|)^{\frac{1}{2}}.
\end{align*}
\end{proof}

{\it Proof of} $(1)\Rightarrow(2)$ {\it of Theorem \ref{main} }.
By the assumption, we can pick up a net $\{x_j\}_{j_\in \mathcal{J}}$ 
in $M_\varphi$ 
which satisfy the conditions of Lemma \ref{lem 9}. 
Now we apply the Lemma \ref{lem 10} to this net for fixed 
$\alpha\in I$, then $\varphi(|X_j(\alpha)_{kl}|)$ converges to 0 
for any $k,l=1,2,\dots,n_\alpha$, where $X_j=\Delta(x_j^2)-1\otimes x_j^2$. 
By Lemma \ref{lem 11}, it implies that 
$\varphi(|X_j(\alpha)|_{kl})$ converges to 0 
for any $k,l=1,2,\dots,n_\alpha$. 
Since we have 
$$(\varphi\otimes\varphi)(|X_j(\alpha)|)
=\sum_{1\leq k,\,l \leq n_\alpha}\varphi(e(\alpha)_{kl})
\varphi(|X_j(\alpha)|_{kl}),$$
we see $(\varphi\otimes\varphi)(|X_j(\alpha)|)$ converges to 0. 
By Lemma \ref{lem 12}, we see 
$\|W^*(\Lambda(e(\alpha)_{k1})\otimes\Lambda(x_j))
          -\Lambda(e(\alpha)_{k1})\otimes\Lambda(x_j)\|$ converges to 0 
for any $k=1,2,\dots,n_\alpha$. 
Then we have 
\begin{align*}
&\big\|W^*(\Lambda(e(\alpha)_{kl})\otimes\Lambda(x_j))
          -\Lambda(e(\alpha)_{kl})\otimes\Lambda(x_j)\big\|\\
=&\,\big\|(J\sigma_{i\over2}^\varphi(e(\alpha)_{1l})^* J\otimes1)
  (W^*(\Lambda(e(\alpha)_{k1})\otimes\Lambda(x_j))
          -\Lambda(e(\alpha)_{k1})\otimes\Lambda(x_j))\big\|\\
\leq&\,\big\|(J\sigma_{i\over2}^\varphi(e(\alpha)_{1l})^* J\otimes1)\big\|
      \big\|W^*(\Lambda(e(\alpha)_{k1})\otimes\Lambda(x_j))
          -\Lambda(e(\alpha)_{k1})\otimes\Lambda(x_j)\big\|. 
\end{align*}
This implies that $\|W^*(\Lambda(e(\alpha)_{kl})\otimes\Lambda(x_j))
          -\Lambda(e(\alpha)_{kl})\otimes\Lambda(x_j)\|$ 
converges to 0 for any $k,l=1,2,\dots,n_\alpha$. 
By taking a linear combination, 
$\|W^* (z_\alpha\eta\otimes\Lambda(x_j))
          -z_\alpha\eta\otimes\Lambda(x_j)\|$ 
converges to 0 for any vector $\eta\in H$. 
Take a vector $\eta\in H$.
For any $\varepsilon>0$, there exists a finite subset $F$ of $I$ such that 
$\|\sum_{\alpha\in F}z_\alpha\eta -\eta\|<\varepsilon$. 
By the above arguments, we can take $j_0$ in $\mathcal{J}$ such that 
$\sum_{\alpha\in F}\|W^* (z_\alpha\eta\otimes\Lambda(x_j))
          -z_\alpha\eta\otimes\Lambda(x_j)\|<\varepsilon$ 
for $j\geq j_0$.
Then we have 
\begin{align*}
&\,\|W^*(\eta\otimes\Lambda(x_j))-\eta\otimes\Lambda(x_j)\|\\
=&\,\Big\|W^*\Big((\eta-\sum_{\alpha\in F}z_\alpha\eta)\otimes\Lambda(x_j)\Big)
           -\Big(\eta-\sum_{\alpha\in F}z_\alpha\eta\Big)\otimes\Lambda(x_j)\\
      &+W^*\Big(\sum_{\alpha\in F}z_\alpha\eta\otimes\Lambda(x_j)\Big)
           -\sum_{\alpha\in F}z_\alpha\eta\otimes\Lambda(x_j)\Big\|\\
\leq&\, 2\Big\|\sum_{\alpha\in F}z_\alpha\eta -\eta\Big\|\\
       &+\sum_{\alpha\in F}\|W^* (z_\alpha\eta\otimes\Lambda(x_j))
          -z_\alpha\eta\otimes\Lambda(x_j)\|\\
<&\, 2\varepsilon+\varepsilon=3\varepsilon,
\end{align*}
for $j\geq j_0$. 
Therefore, $\|W^* \eta\otimes\Lambda(x_j)
          -\eta\otimes\Lambda(x_j)\|$ 
converges to 0 for any vector $\eta\in H$. 
This completes the proof of Theorem \ref{main}. $\hfill\square$

Upon ending this paper, 
we mention a part of Ruan's Theorem 
\cite[Theorem4.5]{Ruan1} as a corollary of 
Theorem \ref{main} for Kac algebras, i.e. 
the invariant weight $\hat{\varphi}$ of $(\hat{M},\hat{\Delta})$ 
is a normal tracial state. 
In this case, $\vph$ is also a trace (\cite{Enock Schwartz2}). 

\begin{cor}\label{cor 4}
Let $\M$ be a discrete Kac algebra. 
Then the following statements are equivalent. 
\begin{enumerate}
   \item It has an invariant mean. 
   \item It is strongly Voiculescu amenable. 
   \item The $C^*$-algebra $\hat{A}$ is nuclear. 
   \item The von Neumann algebra $\hat{M}$ is injective. 
\end{enumerate}
\end{cor}
\begin{proof}
$(1)\Rightarrow(2)\Rightarrow(3)$. 
This has been already proved in Theorem \ref{main}.
$(3)\Rightarrow(4)$. It is trivial. 

$(4)\Rightarrow(1)$. 
Let $E$ be a conditional expectation from $B(H)$ onto $\hat{M}$. 
Note that $\hat{\varphi}$ is a normal trace on $\hat{M}$. 
Take a complete orthonormal system $\{e_p\}_{p\in\mathcal{P}}$. 
Then for any operator $x\in M$ and for any vector $\xi\in H$, we have 
\begin{align*}
\omega_{\xi}* m(x)=&\,\hat{\varphi}\Big(E\Big(\sum_{p\in\mathcal{P}}
                  (\omega_{\xi,\,e_p}\otimes\iota)(W)^*
                   x(\omega_{\xi,\,e_p}\otimes\iota)(W)\Big)\Big)\\
                     =&\,\sum_{p\in\mathcal{P}}\hat{\varphi}
                         \big((\omega_{\xi,\,e_p}\otimes\iota)(W)^* E(x)
                           (\omega_{\xi,\,e_p}\otimes\iota)(W)\big)\\
                     =&\,\sum_{p\in\mathcal{P}}\hat{\varphi}
                         \big((\omega_{\xi,\,e_p}\otimes\iota)(W)
                          (\omega_{\xi,\,e_p}\otimes\iota)(W)^* E(x)\big)\\
                     =&\,\sum_{p\in\mathcal{P}}\hat{\varphi}
                         \big((\omega_{\hat{J}\xi,\,\hat{J}e_p}\otimes\iota)(W)^\ast
                          (\omega_{\hat{J}\xi,\,\hat{J}e_p}\otimes\iota)(W) E(x)\big)\\
                     =&\,\omega_{\hat{J}\xi}(1)\hat{\varphi}(E(x))\\
                     =&\,\omega_{\xi}(1)m(x).
\end{align*}
Therefore, $m$ is a left invariant mean on $M$. 
\end{proof}

As we have seen, 
nuclearity of a compact Kac algebra leads amenability of 
the dual discrete Kac algebra, however, 
it is now open whether it holds in the case of a compact quantum group or not.

\end{document}